\documentclass[a4paper,12pt]{article} % To be LaTeX'ed with cjp3.sty
\def\nn{\nonumber}
\def\slt{{\cal U}_{h,g}(sl(1/2))}

\def\glt{{\cal U}_{h,g}(gl(2))}

\begin{document}

\begin{flushright}
OWUAM-029   \\
July 10, 1998
\end{flushright}

\vspace{2cm}

\renewcommand{\thefootnote}{\fnsymbol{footnote}}

\begin{center}
{\large \bf DRINFELD TWIST FOR TWO-PARAMETRIC DEFORMATION OF 
$ gl(2) $ AND $ sl(1/2) $}\footnote[2]{Presented at the 7th 
International Colloquium "Quantum Groups and Integrable Systems", 
Prague, 18--20 June 1998.}

\vspace{2cm}

N. Aizawa

\medskip
Department of Applied Mathematics, Osaka Women's University

Sakai, Osaka 590-0035, Japan
\end{center}

\vfill

\begin{abstract}
Drinfeld twist is applied to the Lie algebra $gl(2)$ so that 
a two-parametric deformation of it is obtained, which is identical 
to the Jordanian deformation of the $ gl(2) $ obtained by Aneva $et \ al$. 
The same twist element is applied to deform the Lie superalgebra $ sl(1/2) $, 
since the $ gl(2) $ is embedded into the $ sl(1/2) $. By making use of 
the FRT-formalism, we construct a deformation of the Lie supergroup $ SL(1/2). $
\end{abstract}

\newpage
\section{Introduction}

  It is well known that deformation of Lie groups or Lie algebras is not 
unique. As for $ gl(2), $  its multiparametric quantum deformation is 
classified recently \cite{bhp}. In this note, a two-parametric 
deformation of $ gl(2) $ by the Drinfeld twist \cite{dri} is considered 
and it is shown that the 
deformed $ gl(2) $ is identical to the one obtained before \cite{adm}. 
The twist for $ gl(2) $ used here is also applicable to deform the Lie 
superalgebra $ sl(1/2) $, since $ gl(2) \subset sl(1/2). $  
The universal R-matrix for the deformed $ sl(1/2) $ is constructed 
according to the method of twisting. 
This enables us, using the FRT-formalism, 
to construct a deformation of the Lie supergroup $ SL(1/2). $ 
The same scenario 
has been carried out for a one-parametric deformation of $ sl(2) $ and 
$ osp(1/2) $ (Note that $ sl(2) \subset osp(1/2) $)\cite{ck}. The 
present work follows the line adapted in \cite{ck}.

  The deformation by the Drinfeld twist has been developed in recent years. 
For example, multiparametric twists for the Drinfeld-Jimbo deformation 
of simple Lie algebras \cite{resh}, one-parametric twist for $ sl(2) $ 
\cite{ohn}, and twist for Poincar\'e algebra \cite{lrt}, 
Heisenberg algebra \cite{heisen}, esoteric quantum groups \cite{km}, 
$ sl(N) $ \cite{klm}, and $ osp(1/2) $ \cite{ck,super} have been 
considered.

\section{Drinfeld twist}

 This section is devoted to a brief review of the Drinfeld twist. 
Drinfeld develops it in his study of quasi Hopf algebras, we however restrict 
ourselves to ordinary Hopf algebras throughout this note.

 Let $ {\cal A} $ be a Hopf algebra with coproduct $ \Delta_0 $, counit 
$ \epsilon_0 $ and antipode $ S_0. $ Let $ {\cal F} $ be an invertible 
element in $ {\cal A} \otimes {\cal A} $ satisfying the conditions
\begin{eqnarray}
 & &  (\epsilon_0 \otimes id )({\cal F})  = (id \otimes \epsilon_0)({\cal F})
        = 1, \nn \\
 & &  {\cal F}_{12}(\Delta_0 \otimes id)({\cal F}) = {\cal F}_{23}(id \otimes 
        \Delta_0)({\cal F}), 
\end{eqnarray}
then we obtain a new Hopf algebra $ {\cal H} $ with the same algebraic 
relations as $ {\cal A}$. However, Hopf algebraic mappings of $ {\cal H} $ 
are different, they are twisted by the twist element $ {\cal F} $. The 
coproduct, the counit and the antipode for $ {\cal H} $ are given by
\begin{equation}
 \Delta = {\cal F} \Delta_0 {\cal F}^{-1}, \quad
 \epsilon = \epsilon_0, \quad
 S = u S_0 u^{-1},                         \label{hopf}
\end{equation}
where $ u = m(id \otimes S_0)({\cal F}). $

  When the algebra $ {\cal A} $ has a universal R-matrix $ {\cal R}_0 $, the 
universal R-matrix $ {\cal R} $ for $ {\cal H} $ is given by
\begin{equation}
 {\cal R} = {\cal F}_{21} {\cal R}_0 {\cal F}^{-1}. \label{unir}
\end{equation}
This shows that if $ {\cal A} $ is cocommutative, then $ {\cal H} $ is 
a triangular Hopf algebra. This is the case that we would like to 
consider, since we shall start with a Lie algebra.

\setcounter{equation}{0}
\section{Twist for $ gl(2) $}

  The Lie algebra $ gl(2) $ has elements $ Z, H, $ and $ X_{\pm} $ and 
defined by the relations
\begin{equation}
  [H,\; X_{\pm}] = \pm 2X_{\pm}, \quad
  [X_+,\; X_-] = H, \quad
  [Z,\; \bullet] = 0.            \label{relgl2}
\end{equation}
A twist element for $ gl(2) $ is given by
\begin{equation}
  {\cal F}  = \exp\left(\frac{g}{2h} \sigma \otimes Z \right) 
       \exp\left(-\frac{1}{2} H \otimes \sigma \right)_,   \label{tw}
\end{equation}
with
\[
  \sigma \equiv -\ln(1-2hX_+),
\]
where $ h $ and $ g $ are deformation parameters. 

 Acoording to (\ref{hopf}), the coproduct and the antipode are deformed and 
they read
\begin{eqnarray}
 & & \Delta(Z) = Z \otimes 1 + 1 \otimes Z, \nn \\
 & & \Delta(H) = H \otimes e^{\sigma} + 1 \otimes H + {\textstyle \frac{g}{h}} 
     (1-e^{\sigma}) \otimes Z e^{\sigma}, \nn \\
 & & \Delta(X_+) = X_+ \otimes 1 + e^{-\sigma} \otimes X_+, \nn \\
 & & \Delta(X_-) = X_- \otimes e^{\sigma} + 1 \otimes X_- 
     - hH \otimes e^{\sigma} H 
 - {\textstyle \frac{h}{2}} H (H+2) \otimes e^{\sigma}(e^{\sigma}-1) \nn \\
 & & \hspace{1.2cm}    + g(e^{\sigma}-1)\otimes Ze^{\sigma} H 
 + g(H-e^{\sigma} + 1) \otimes Z e^{\sigma} \nn \\ 
 & & \hspace{1.2cm} + g(e^{\sigma}-1)(H+e^{\sigma}+1) \otimes Z e^{2\sigma} \nn \\
 & & \hspace{1.2cm} - {\textstyle \frac{g^2}{2h}}(e^{\sigma}-1) \otimes Z^2 e^{\sigma}
     - {\textstyle \frac{g^2}{2h}} (e^{\sigma}-1)^2 \otimes Z^2 e^{2\sigma}, \nn \\
  & & S(Z) = -Z, \nn \\
  & & S(H) = -H e^{-\sigma} + {\textstyle \frac{g}{h}} Z(e^{-\sigma}-1), \nn \\
  & & S(X_+) = - X e^{\sigma}, \nn \\
  & & S(X_-) = -\left\{ X_- + {\textstyle \frac{h}{2}} H^2 (e^{-\sigma}+1) 
  - hH(e^{-\sigma}-1) - gHZe^{-\sigma} \right. \nn \\
  & & \hspace{1.1cm} \left.  + gZ(e^{-\sigma}-1) 
      + {\textstyle \frac{g^2}{2h}} (e^{-\sigma}-1)Z^2 \right\}e^{-\sigma}. 
      \label{hopfgl2}
\end{eqnarray}
The universal R-matrix for the non-cocommutative coproduct given above reads
\begin{equation}
{\cal R}  =  \exp\left(\frac{g}{2h}Z \otimes \sigma\right) 
     \exp\left(-\frac{1}{2}\sigma \otimes H \right) 
     \exp\left(\frac{1}{2}H \otimes \sigma\right)
     \exp\left(-\frac{g}{2h}\sigma \otimes Z \right)_.       \label{unirgl2}
\end{equation}
A counit is necessary for a Hopf algebra, it is undeformed, that 
is, $ \epsilon(Z)=\epsilon(H)=\epsilon(X_{\pm}) = 0. $ 
Thus we arrive to the definition

\medskip
\noindent
{\bf Definition 1.} The triangular Hopf algebra generated by 
$ \{ Z, H, X_{\pm} \} $ satisfying the relations (\ref{relgl2}) 
and (\ref{hopfgl2}) is said to be the two-parametric deformation of 
$ {\cal U}(gl(2)) $ by twisting or $ \glt. $

\medskip
  Let us take the particular nonliear combinations of generators,
\begin{eqnarray}
 A &=& Z, \nn \\
 H' &=& e^{-\sigma /2} H \nn \\
 X &=& \frac{1}{2h} \sigma, \label{nonlin} \\
 Y &=& e^{-\sigma /2}\left(X_- + \frac{h}{2} H^2 \right) 
      - \frac{h}{8}e^{\sigma /2} (e^{-\sigma}-1). \nn
\end{eqnarray}
Then $ A, H', X $ and $ Y $ satisfy the commutation relations,
\begin{eqnarray}
& & [X,\; Y] = H', \qquad
      [H', \; X] = 2{\sinh hX \over h},
      \nn \\
 & &  [H', \; Y] = -Y (\cosh hX) - (\cosh hX) Y,  \label{jorcomm} \\
 & & [A, \; \bullet] = 0. \nn
\end{eqnarray}
The Hopf algebra mappings for these generators are given by
\begin{eqnarray}
& & \Delta(A) = A \otimes 1 + 1 \otimes A, \nn \\
& & \Delta(H') = H' \otimes e^{hX} + e^{-hX} \otimes H'
    - \frac{2g}{h} \sinh hX \otimes A e^{hX}, \nn \\
& & \Delta(X) = X \otimes 1 + 1 \otimes X, \nn \\
& & \Delta(Y) = Y \otimes e^{hX} + e^{-hX} \otimes Y 
- \frac{g^2}{h} \sinh hX \otimes A^2 e^{hX} + gH' \otimes Ae^{hX}, \nn \\
& &   \epsilon(\bullet) = 0, \nn \\
 & & S(A) = -A, \qquad S(X) = -X, \nn \\
 & & S(H') = -e^{hX}H'e^{-hX} - \frac{2g}{h} (\sinh hX) A, \nn \\
 & & S(Y) = -e^{hX} Y e^{-hX} + \frac{g^2}{h}(\sinh hX) A^2 
 + g e^{hX}H'Ae^{-hX}. \nn
\end{eqnarray}
Therefore the algebra generated by $ \{A, H', X, Y \} $ is nothing but the 
one introduced in \cite{adm}. The authors of \cite{adm} define the algebra 
so as to be dual to the Jordanian matrix quantum group $ GL_{h,g}(2) $
\cite{dmmz,ag}.

\setcounter{equation}{0}
\section{Twist for $sl(1/2)$}

  The Lie superalgebra $ sl(1/2) $ has four even and four odd elements denoted by 
$ Z, H, X_{\pm} $ and $ v_{\pm}, \bar v_{\pm} $, respectively. They satisfy 
the relations
\begin{eqnarray}
 & & [H,\; X_{\pm}] = \pm 2 X_{\pm}, \quad \ \,
     [X_+,\; X_-] = H, \qquad [Z,\; \bullet] = 0, \nn \\
 & & [X_{\pm},\; v_{\mp}] = -v_{\pm}, \qquad
 [X_{\pm},\; \bar v_{\mp}] = \bar v_{\pm},  \nn \\ 
 & &  
 [Z, \; v_{\pm}] = v_{\pm},\qquad \hspace{6mm}  [Z,\; \bar v_{\pm}] = -\bar v_{\pm}, \nn \\
 & & [H, \; v_{\pm}] = \pm v_{\pm}, \qquad \ \,
 [H, \; \bar v_{\pm}] = \pm \bar v_{\pm}, \qquad 
 [X_{\pm},\; v_{\pm} ] = [X_{\pm},\; \bar v_{\pm}] = 0, \nn  \\
 & & \{ v_{\pm},\; v_{\pm} \} = \{ \bar v_{\pm},\; \bar v_{\pm} \} = 
     \{ v_{\pm},\; v_{\mp} \} = \{ \bar v_{\pm},\; \bar v_{\mp} \} = 0, \nn \\ 
 & & \{ v_{\pm}, \; \bar v_{\pm} \} = X_{\pm}, \qquad 
 \{ \bar v_{\pm},\; v_{\mp} \} = {\textstyle 1 \over 2} (Z \pm H). \label{relsl}
\end{eqnarray}
  It is easily seen from the above relations that the 
even elements $ \{Z, H, X_{\pm} \} $ form a $ gl(2) $ subalgebra and 
the universal enveloping algebra $ {\cal U}(sl(1/2)) $ is generated by 
$ \{ Z, H, X_{\pm}, v_{-}, \bar v_{+} \} $. The observation of 
$ gl(2) \subset sl(1/2) $ implies that the twist element for $ gl(2) $ 
can be used to twist $ sl(1/2) $. 
Using (\ref{tw}), the twisted coproduct for the odd elements are 
given by
\begin{eqnarray}
  \Delta(\bar v_+) &=& \bar v_+ \otimes e^{-\sigma/2} 
      + \exp\left(-{\textstyle \frac{g}{2h}} \sigma \right) \otimes \bar v_+, \nn \\
  \Delta(v_-) &=& v_- \otimes e^{\sigma/2} 
      + \exp\left( {\textstyle \frac{g}{2h}}\sigma\right) \otimes v_-  
      + hH \exp\left({\textstyle \frac{g}{2h}}\sigma \right) \otimes v_+e^{\sigma} \nn \\  
      &-& g v_+ e^{\sigma} \otimes Z e^{\sigma/2}
  - g (e^{\sigma}-1) \exp\left({\textstyle \frac{g}{2h}}\sigma \right) \otimes Z v_+ e^{\sigma}, 
  \label{copsl}
\end{eqnarray}
and the antipode is
\begin{eqnarray}
  & & S(\bar v_+) = - \bar v_
      + \exp\left({\textstyle \frac{1}{2h}} (h+g) \sigma\right) \nn \\
  & & S(v_-) = - (v_- - hHv_+ + gv_+(1+Z))
      \exp\left(-{\textstyle \frac{1}{2h}}(h+g) \sigma \right). \label{apsl}
\end{eqnarray}
The counit is undeformed and given by 
$ \epsilon(\bar v_+) = \epsilon(v_-) = 0. $ 
The Hopf algebra mappings for even elements have already been given in 
(\ref{hopfgl2}). 

\medskip
\noindent
{\bf Definition 2.} The triangular Hopf algebra generated by 
$ \{ Z, H, X_{\pm}, \bar v_+, v_- \} $ satisfying the relations (\ref{relsl}), 
(\ref{copsl}) and (\ref{apsl}) is said to be the two-parametric deformation of 
$ {\cal U}(sl(1/2)) $ by twisting or $ \slt. $

\medskip
The universal R-matrix for $ \slt $ is same as the one for $ \glt. $ 
Noting that the fundamental representation of $ \slt $ is same as 
$ sl(1/2) $, we obtain the R-matrix in the fundamental representation 
of $ \slt $, and see that it is a direct sum of four matrices
\begin{equation}
  R = (1) \oplus \check R \oplus \check R^{-1} \oplus \bar R, 
  \label{rsl2}
\end{equation}
where
\begin{equation}
  \check R = \left(
\begin{array}{rr}
  1 & 2g  \\
  0 & 1   \\
\end{array}
\right),
\qquad
  \bar R = \left(
\begin{array}{rrrr}
  1 & h+g & -h-g & h^2 - g^2  \\
  0 & 1   &  0   & h-g        \\
  0 & 0   &  1   & -h+g       \\
  0 & 0   &  0   &  1
\end{array}
\right)_.
\end{equation}
The matrix $ \bar R $ is the R-matrix (\ref{unirgl2}) in the fundamental 
representation of $ gl(2). $

\setcounter{equation}{0}
\section{Two-parametric deformation of $ SL(1/2) $}

  Using the R-matrix (\ref{rsl2}) and $ {\bf Z}_2 $ graded version of 
FRT-formalism \cite{FRT}, a matrix quantum supergroup dual to $ \slt $ 
can be constructed. 

Introducing a $ 3 \times 3 $ supermatrix 
\[
  M = \left(
    \begin{array}{ll}
    e & \Psi \\ \Theta & T
    \end{array}
  \right)_,
\]
where
\[
  \Psi = (\xi, \eta), \qquad  
  \Theta = \left(
   \begin{array}{c}
     \gamma \\ \delta
   \end{array}
  \right)_,
  \qquad 
  T = 
  \left(
  \begin{array}{cc}
   a & b \\ c & d 
  \end{array}
  \right)_,
\]
and $ e, T $ are even elements and $ \Psi, \Theta $ are odd elements, 
the FRT-formalism guarantees that the commutation relations for 
the entries of $ M $ are given by RMM-relation and their Hopf algebra 
mappings are given by
\begin{equation}
    \Delta(M) = M \otimes M,\quad
  \epsilon(M) = I_3, \quad
  S(M) = M^{-1}. \label{hopfg}
\end{equation}
The commutation relations for the entries of $ M $ read
\begin{equation}
  \begin{array}{lll}
  e \Psi = \Psi e \check R, 
  &  \quad 
  \check R e \Theta = \Theta e,  
  & \quad 
  \check R e T = T e \check R, 
  \\ 
  ( \xi \Psi \ \eta \Psi ) = - (\Psi \xi \ \Psi \eta ) \bar R, 
  & \quad 
  \bar R \left(
  \begin{array}{c}
   \gamma \Theta \\ \delta \Theta
  \end{array}\right)
  =
  - \left(
  \begin{array}{c}
   \Theta \gamma \\ \Theta \delta
  \end{array}
  \right)_,
  & 
  \\ 
  \check R (\xi \Theta \ \eta \Theta ) \check R 
  =
  -\left(
   \begin{array}{c}
     \gamma \Psi \\ \delta \Psi
   \end{array}\right)_, 
  & \quad 
   \bar R \left(
     \begin{array}{c}
      \gamma T \\ \delta T
     \end{array}
   \right)
   =
   \left(
    \begin{array}{c}
     T \gamma \\ T \delta
    \end{array}
   \right)
   \check R, 
   & 
   \\
   \check R ( \xi T \ \eta T ) 
   =
   (T \xi \ T \eta) \bar R,
   & \quad 
   \bar R T_1 T_2 = T_2 T_1 \bar R.  & 
 \end{array}
 \label{relg}
\end{equation}
The last relation shows that the submatrix $ T $ satisfy the same 
algebraic relations as the $ GL_{h,g}(2) $ \cite{ag}. As in \cite{ag}, 
a determinant for $ T $ is defined by $ det T = ad - bc - (h+g)ac $, 
then it is not difficult to see that $ det T $ is not a center of 
deformed $ SL(1/2) $ and the noncommutativity is independent of $h$
\begin{eqnarray}
 & & [\Omega, \; detT] = 0, \quad \Omega = e, \; \xi,\; \delta,\; c \nn \\
 & & [\eta, \; detT] = -2g \, \xi\, detT, \quad \quad 
     [\gamma, \; detT] = 2g \,\delta\, detT, \nn \\
 & & [a,\; detT] = [detT,\; d] = 2g\, c\, detT, \nn \\
 & & [b,\; detT] = 2g \{ (detT) d - a (detT) \}. \nn
\end{eqnarray}
Assuming that the $ det T $ is invertible, the explicit form of the 
inverse matrix of $ T $ can be obtain (see \cite{ag} for the formula). 

  We define a superdeterminant for the quantum supermatrix $ M$ by
\begin{equation}
 sdet M = (detT)^{-1} (e - \Psi T^{-1} \Theta). \label{sdet}
\end{equation}
This has the same form as the undeformed case (it is also called 
Berezinian in undeformed case). Direct computation shows that 
the $ sdet M $ commute with all elements of the deformed $ SL(1/2) $ 
so that we can safely set $ sdet M = 1. $ A coproduct and a counit 
for $ M $ are obvious from (\ref{hopfg}), however, an antipode is 
not. It is necessary to assume that the combination 
$ e - \Psi T^{-1} \Theta $ has a inverse. Then the inverse matrix of $ M$ is 
given by \cite{ck}
\begin{equation}
  M^{-1} = 
  \left(
    \begin{array}{cc}
      1 & 0 \\ -T^{-1} \Theta & I_2
    \end{array}
  \right)
  \left(
    \begin{array}{cc}
      (e - \Psi T^{-1} \Theta)^{-1} & 0 
      \\
      0 & T^{-1}
    \end{array}
  \right)
   \left(
    \begin{array}{cc}
      1 & -\Psi T^{-1} 
      \\
      0 & I_2
    \end{array}
  \right)_,
  \label{inverse}
\end{equation}
where $ I_2 $ is the $ 2 \times 2 $ unit matrix. 

\medskip
{\bf Definition 3.} An algebra generated by the entries of $ M $ satisfying 
(\ref{relg}),  (\ref{hopfg}), (\ref{inverse}) and $ sdet M = 1 $ 
%with $ sdet M = 1 $ 
is said to be the two-parametric deformation of  $ Fun(SL(1/2)) $ or 
$ SL_{h,g}(1/2). $

\bigskip

 We finally give some remarks. The twist element $ {\cal F} $ (\ref{tw}) 
is also applicable to deform $ osp(1/2) \oplus u(1) $, since 
$ gl(2) \subset osp(1/2) \oplus u(1) $. The universal R-matrix for the 
obtained algebra is used to deform the supergroup $ OSp(1/2) \otimes U(1) $. 
The inclusion of a odd element of $ sl(1/2) $ into a twist element may 
be possible. The inclusion of a odd elements is found for $ osp(1/2) $ 
recently \cite{super}.

\end{document}